\documentclass[11pt,a4paper]{article}
\usepackage[pagebackref=true,linktoc=page,colorlinks=true,linkcolor=blue,citecolor=blue]{hyperref}
\usepackage{setspace}
\usepackage{epsfig}
\usepackage{psfrag}
\usepackage[retainorgcmds]{IEEEtrantools} 
\usepackage{fullpage,color}
\usepackage{authblk}
\usepackage{latexsym}
\usepackage{amsmath,amssymb,amsfonts,amsthm}
\usepackage{graphicx}
\usepackage{natbib}
\usepackage{pstool,psfrag}
\usepackage[section]{placeins}
\newtheorem{theorem}{Theorem}
\newtheorem{lemma}{Lemma}

\makeatletter
 
\newcommand{\Xf}{\mbox{$X_0$}} 
\newcommand{\Xs}{\mbox{$X_1$}} 
\newcommand{\Xt}{\mbox{$X_2$}} 
\newcommand{\indf}{0} 
\newcommand{\inds}{1} 
\newcommand{\indt}{2} 
\newcommand{\locs}{\mbox{$\beta_1$}} 
\newcommand{\scas}{\mbox{$\alpha_1$}} 
\newcommand{\loct}{\mbox{$\beta_2$}} 
\newcommand{\scat}{\mbox{$\alpha_2$}}

\definecolor{light-gray}{gray}{0.80} 
\title{\textbf{Conditional independence and conditioned limit laws}}

\author[$$]{\textbf{Ioannis Papastathopoulos}} \affil[$$]{\small School
  of Mathematics and Maxwell Institute, University of Edinburgh,
  Edinburgh, EH9 3FD} \affil[$$]{\small i.papastathopoulos@ed.ac.uk}
 \date{}
 \date{}
\begin{document}

\maketitle
% \tableofcontents
\begin{abstract}
  Conditioned limit laws constitute an important and well developed
  framework of extreme value theory that describe a broad range of
  extremal dependence forms including asymptotic independence.\ We
  explore the assumption of conditional independence of $X_{\inds}$
  and $X_{\indt}$ given $X_{\indf}$ and study its implication in the
  limiting distribution of $(X_{\inds},X_{\indt})$ conditionally on
  $X_{\indf}$ being large.\ We show that under random norming,
  conditional independence is always preserved in the conditioned
  limit law but might fail to do so when the normalisation does not
  include the precise value of the random variable in the conditioning
  event.\
\end{abstract}

% \section*{Notes}
% \begin{itemize}
% \item notes for paper: Resnick's joint limit is a mixing
%   distribution of HT marginal distributions.
% \item So to estimate we need to have HT marginals distributions
%   whereas the Resnick's marginals will not help (?)
% \item Think of inverse transform of link between Resnick and HT
%   marginal distributions.
% \item If they can be inverted the limit joint distribution is an
%   integral so it is more complicated. Resnick joint distribution is
%   not a product of Resnick marginal distribution.
% \item Mathematically more elegant to use Resnick approach for the
%   norming constants.
% \end{itemize}

\section{Introduction}
\label{sec:intro}
Extreme value theory is a highly active area of research and its
methods and applications are the epitome of risk modelling and
statistical estimation of rare events.\ Extreme events that occur in a
broad range of disciplines such as in environmental processes or in
finance and insurance are typically multivariate in nature and usually
have tremendous socio-economic impact.\ The recent technological
advances have resulted in an ever-increasing amount of information
available across the whole spectrum of applied sciences.\ As such,
when modelling data in several dimensions, one is typically confronted
by the \emph{curse of dimensionality}.\ It is widely recognized that
the construction of more efficient statistical models and techniques
that overcome this problem is imperative.\ \emph{Conditional
  independence} constitutes one of the most fundamental tools and
concepts in this direction
\citep{besa74,dawi79,JoeWhittaker,lau96,coxwerm96}.\ On the other
hand, the central concept of regular variation, its extensions and
refinements, provide the recipe for the development of asymptotically
justified extreme value models.\ The purpose of this short note is to
illustrate some implications of conditional independence in
conditioned limit laws \citep{hefftawn04,heffres07}, a key and
well-developed framework that embodies a broad range of extremal
dependence forms.
 
For ease of exposition, consider a random vector $(\Xf,\Xs,\Xt)$ in
$\mathbb{R}^{3}$.\
% \marginpar{\calt{On right lines but could be strengthened by improving
%     text flow}}
The result developed in this paper extends to higher-dimensional
settings for $\Xs$ and $\Xt$ in a straightforward manner so this is
not restrictive.\ We assume that $\Xs$ and $\Xt$ are conditionally
independent given $\Xf$.\ Informally, this means that the conditional
distribution of $\Xs$ given $(\Xf,\Xt)$ is equal to the conditional
distribution of $\Xs$ given $\Xf$ alone.\ In other words, once $\Xf$
is known any further information about $\Xt$ is irrelevant to
uncertainty about $\Xs$ and therefore it readily follows that for any
$x_{\inds},x_{\indt} \in \mathbb{R}$,
\begin{equation}
  \mathbb{P}\left(\Xs<x_{\inds},\Xt<x_{\indt}\mid \Xf\right)= \mathbb{P}\left(\Xs<x_{\inds}\mid
    \Xf\right) \mathbb{P}\left(\Xt<x_{\indt}\mid \Xf\right),
  \label{eq:CI}
\end{equation}
almost surely.\ 
% Multivariate extreme value theory focuses upon asymptotic
% characterisations of joint tail region of random vectors.\ Let
% $F_X,F_Y$ and $F_Z$ denote the distribution function of $X,Y$ and
% $Z$, respectively, and define the coefficient of \emph{asymptotic
% dependence}
We show that conditional independence may naturally be
\emph{preserved} in limiting laws of random vectors given an extreme
component \citep{hefftawn04,heffres07}.\ Conditioned limit laws
provide a rich description of extremes of random vectors that do not
necessarily grow at the same rate and may exhibit \emph{asymptotic
  independence} which means that the coefficient
\begin{equation}
  \chi=\lim_{p\rightarrow
    1-}\mathbb{P}\left[F_{\scriptsize \Xs}(\Xs)>p,F_{\scriptsize \Xt}(\Xt)>p\mid
    F_{\scriptsize \Xf}(\Xf)> p\right],
\end{equation}
can be 0.\ % The cases $\chi>0$ and $\chi=0$ correspond to asymptotic
% dependence and asymptotic independence of $(X,Y,Z)$.\
Conditioned limit laws were systematically studied for first time by
\cite{hefftawn04} who examined the limiting conditional distribution
of affinely transformed random vectors as the conditioning variable
becomes large.\ Assuming identical margins with distribution function
being asymptotically equivalent to the unit exponential distribution,
i.e., $F_{\scriptsize \Xf}(x)=F_{\scriptsize \Xs}(x)=F_{\scriptsize
  \Xt}(x)\sim 1-\exp(-x)$, as $x\rightarrow \infty$, Heffernan and
Tawn showed that for a broad range of dependence structures of
$(\Xf,\Xs,\Xt)$, there exist scaling functions
$\scas,\scat:(0,\infty)\mapsto (0,\infty)$, location functions
$\locs,\loct:(0,\infty)\mapsto \mathbb{R}$ and a joint distribution
$G$ on $[-\infty,\infty)\times [-\infty,\infty)$ with non-degenerate
marginals, such that as $t\rightarrow \infty$
% and subject to existence of
% densities
\begin{equation}
  \mathbb{P}\left[\frac{\Xs-\locs(\Xf)}{\scas(\Xf)}<x_{\inds},
    \frac{\Xt-\loct(\Xf)}{\scat(\Xf)}<x_{\indt}~\Big|~\Xf>t\right]
  \overset{\mathcal{D}}{\rightarrow} G(x_{\inds},x_{\indt}),
  \label{eq:CEVRN}
\end{equation}
on $[-\infty,\infty]\times[-\infty,\infty]$ where
$\overset{\mathcal{D}}{\rightarrow}$ stands for weak
convergence. % where $G$ is a joint
% distribution function with non-degenerate marginals that place no
% mass at $+\infty$, i.e., $G_{\inds}(x) = G_{\indt}(x)=1$, as
% $x\rightarrow \infty$, and $G_{\inds}(y)=G(y,\infty)$ and
% $G_{\indt}(z)=G(\infty,z)$.\
Although the original formulation of Heffernan and Tawn relied on the
existence of densities, formulation~\eqref{eq:CEVRN} is presented here
in the compact form of \cite{heffres07} who provided a formal
mathematical examination of the more general conditioned limit
formulation that there exists a joint distribution $H$ on
$[-\infty,\infty)\times [-\infty,\infty)$ with non-degenerate
marginals such that as $t\rightarrow \infty$
\begin{equation}  
  \mathbb{P}\left[\frac{\Xs-\locs(t)}{\scas(t)}<x_{\inds},\frac{\Xt-\loct(t)}{\scat(t)}<x_{\indt}~\Big|~
    \Xf>t\right] \overset{\mathcal{D}}{\rightarrow} H(x_{\inds},x_{\indt}),
  \label{eq:CEV}
\end{equation} 
on $[-\infty,\infty]\times[-\infty,\infty]$, subject to the sole
assumption of $\Xf$ belonging to the domain of attraction of an
extreme value distribution.\ Expressions~\eqref{eq:CEVRN} and
\eqref{eq:CEV} can be rephrased more generally as special cases of
joint probability convergence; here we use the conditional
representation to highlight the connection with conditional
independence.

Limit expressions~\eqref{eq:CEVRN} and \eqref{eq:CEV} differ in the
way $\Xs$ and $\Xt$ are normalised since in
expression~\eqref{eq:CEVRN}, the precise value of $\Xf$ that occurs
with $\Xf>t$ is used, whereas in expression~\eqref{eq:CEV} only
partial information about the value of $\Xf$ is exploited since only
the level value $t$ that $\Xf$ exceeds is used.\ Following the
terminology of \cite{heffres07} we refer to limit
expressions~\eqref{eq:CEVRN} and \eqref{eq:CEV} as the conditional
extreme value model with \emph{random norming} and \emph{deterministic
  norming}, respectively, and for the remaining part we assume without
loss of generality that $\Xf$ has a unit Pareto distribution, i.e.,
$\mathbb{P}(\Xf<x) = 1-1/x$, $x>1$.\ The unit Pareto marginal scale
for $\Xf$ is primarily chosen for convenience but the result of this
paper can be stated, with modified proofs (\cite{kulisoul}), for $\Xf$
being in the domain of attraction of an extreme value distribution.\
We also assume that the two limit expressions~\eqref{eq:CEVRN} and
\eqref{eq:CEV} simultaneously hold for the same pair of norming
functions. In general, this may not always be true but a necessary and
sufficient condition is to assume that $(\alpha_i,\beta_i)$, $i=1,2$,
in expression~\eqref{eq:CEV} are extended regularly varying, i.e.,
there exist $\rho_{\scriptsize \inds},\rho_{\scriptsize
  \indt},\kappa_{\scriptsize \inds},\kappa_{\scriptsize
  \indt}\in\mathbb{R}$ such that as $t\rightarrow \infty$,
\[
\frac{\alpha_i(t x)}{\alpha_i(t)} \rightarrow x^{\rho_i} \qquad
\text{and} \qquad \frac{\beta_i(tx) - \beta_i(t)}{\alpha_i(t)}
\rightarrow \psi_i (x),
\]
for $x>0$, where
\begin{equation}
\psi_i(x) =
\begin{cases}
  \kappa_i \left(x^{\rho_i}-1\right)/\rho_i & \rho_i \neq 0\\
  \kappa_i \log x &\rho=0,
\end{cases}
\label{eq:psi}
\end{equation}
for $i =1,2$ \citep{resnzebe14}.\ % Assuming in addition, and without
% loss of generality, that $F_{\scriptsize \Xf}(x) = 1-1/x$, $x>1$,
  
In Section~\ref{sec:main_thm} we state the main theorem of the paper
which shows that the conditional independence property~\eqref{eq:CI}
is preserved in the conditional extreme value model with random
norming, meaning that for any $x_{\inds},x_{\indt}\in \mathbb{R}$,
\begin{equation}
  G(x_{\inds},x_{\indt})=G_{\inds}(x_{\inds})\,G_{\indt}(x_{\indt}),
  \label{eq:FACTOR}
\end{equation}
where $G_1(x_1)=\lim_{x_2\rightarrow \infty}G(x_1,x_2)$ and
$G_2(x_2)=\lim_{x_1\rightarrow\infty}G(x_1,x_2)$, but might fail to do
so in the conditional extreme value model with deterministic norming.\
In Section~\ref{sec:discussion} we discuss practical consequences of
conditional independence in conditioned limit laws.\ A proof is given
in Section~\ref{sec:proof}.
% \cite{dasresn11a,dasresn11b} \cite{resn08} \cite{kulisoul}
% \cite{janssege14} \cite{resnzebe13,resnzebe14} \cite{fougsoul12}

\section{Main result}
\label{sec:main_thm}
Let $\pi_{\scriptsize \inds},\pi_{\scriptsize \indt}:(0,\infty)\times
\mathcal{B}(\mathbb{\bar{R}})\mapsto [0,1]$ be Markov kernels defined
by
\[
\pi_{\inds}(x,A) = \mathbb{P}(\Xs\in A\mid \Xf=x) \quad \text{and}
\quad \pi_{\indt}(x,A) = \mathbb{P} (\Xt\in A\mid \Xf=x).
\]
Let also $\sigma(\Xf)=\{\Xf^{-1}(B):B\in
\mathcal{B}(\mathbb{R})\}$ be the sigma algebra
generated by $\Xf$.
\begin{theorem}
  \label{thm:main_thm}
  Let $(\Xf,\Xs,\Xt)$ be a random vector from a probability space to
  $\mathbb{R}^3$ and suppose that given $\Xf$, $\Xs$ is conditionally
  independent of $\Xt$, i.e.,
  \[
  \mathbb{E}\left\{f_{\inds}\left(\Xs\right)\,f_{\indt}\left(\Xt\right)\mid
    \sigma(\Xf)\right\} = \mathbb{E}\left\{f_{\inds}\left(\Xs\right)\mid
    \sigma(\Xf) \right\} \mathbb{E}\left\{f_{\indt}\left(\Xt\right)\mid
    \sigma(\Xf) \right\}
  \]
  almost surely, for any Borel measurable functions
  $f_{\inds},f_{\indt}:\mathbb{R}\mapsto\mathbb{R}$ for which
  $f_{\inds}\circ \Xs$ and $f_{\indt} \circ \Xt$ are Lebesgue
  integrable with respect to $\mathbb{P}$.\ Suppose there exist two
  pairs of scaling and location functions
  $(\scas,\locs),(\scat,\loct):(0,\infty)^2\mapsto (0,\infty)\times
  \mathbb{R}$ that are extended regularly varying
  % with real parameters $(\rho_Y,\kappa_Y)$ and $(\rho_Z,\kappa_Z)$,
  % respectively
  and families of non-degenerate probability distribution functions
  $\{G_{v;\tiny \inds},v>0\}$ and $\{G_{v;\tiny \indt},v>0\}$ on
  $[-\infty,\infty]$, such that as $t\rightarrow \infty$
  \begin{equation}
    \pi_{\scriptsize \inds}(t,\left(-\infty,\scas(t) x_{\inds} + \locs(t)\right]) \overset{\mathcal{D}}{\rightarrow}
    G_{\tiny \inds} (x_{\inds}) \quad \text{and} \quad
    \pi_{\scriptsize \indt}(t,\left(-\infty,\scat(t)x_{\indt} + \loct(t)\right]) \overset{\mathcal{D}}{\rightarrow}G_{{\tiny \indt}}
    (x_\indt),
    \label{eq:kernel_conv}
  \end{equation}
  on $[-\infty,\infty]$ where $G_1=G_{1;1}$ and $G_2=G_{1;2}$.\ Then,
  as $t\rightarrow \infty $ and for $f_{\inds}$ and $f_{\indt}$
  bounded
  \begin{IEEEeqnarray}{rCl}
    &&\mathbb{E}\left\{f_{\inds}\left(\frac{\Xs-\locs(\Xf)}{\scas(\Xf)}\right)\,
      f_{\indt}\left(\frac{\Xt-\loct(\Xf)}{\scat(\Xf)}\right)~\Big|~
      \Xf > t\right\} \rightarrow
    \mathbb{E}_{G_{\inds}}\{f_{\inds}(\Xs)\}
    \mathbb{E}_{G_{\indt}}\{f_{\indt}(\Xt)\}.
    \label{eq:week2}
  \end{IEEEeqnarray}  
  and
  \begin{IEEEeqnarray}{rCl}
    &&\mathbb{E}\left\{f_{\inds}\left(\frac{\Xs-\locs(t)}{\scas(t)}\right)\,
      f_{\indt}\left(\frac{\Xt-\loct(t)}{\scat(t)}\right)~\Big|~ \Xf >
      t\right\} \rightarrow \int_{\scriptsize \inds}^{\infty}
    \mathbb{E}_{G_{v;{\tiny \inds}}}\{f_{\inds}(\Xs)\} \mathbb{E}_{G_{v;\tiny
        \indt }}\{f_{\indt}(\Xt)\} v^{-2} ~dv
    \nonumber\\ \label{eq:week1}
  \end{IEEEeqnarray}

  % where $\{G_{v;Y}:0<v<\infty\}$ and $\{G_{v;\Xt}:0<v<\infty\}$ are
  % families
  % of non-degenerate probability distributions on $[-\infty,\infty)$
  % such that $G_{1;Y}=G_{Y}$ and $G_{1;\Xt}=G_{\Xt}$
  % Under the same conditions of Theorem~\ref{thm:main_thm}
\end{theorem}

% \marginpar{ {\ccom Discussion here will be based around the
% statistical implications of the following equation
% and~\eqref{eq:FACTOR} }}

\section{Discussion}
\label{sec:discussion}
Conditional independence shows that important simplifications can be
achieved at a practical level.\ To elaborate on this, consider a
regression setting where interest lies in estimating the conditional
probability on the left hand side of equation~\eqref{eq:CI}.\ A
statistician has virtually two options.\ The first is to estimate the
joint conditional probability by means of either a parametric model or
to proceed non-parametrically.\ The second approach consists of
estimating in the same spirit the lower dimensional marginal
conditional distributions and combine them multiplicatively to yield
an estimate.\ Thus, in light of extra information regarding
conditional independence, the second approach is more efficient since
it relates to a lower dimensional problem.\

Theorem~\ref{thm:main_thm} shows that the inclusion of the precise
value of the random variable in the conditioning event adds enough
detail to the normalisation to always allow the limit law to factorise
under the presence of conditional independence.\ This would offer
substantial simplifications in practical as well as theoretical
extreme value analyses.\ On the other hand, it follows from
Theorem~\ref{thm:main_thm} and Lemma~\ref{lemma:resnzebe} that
\begin{equation}
  H(x_{\inds},x_{\indt}) = \int_{\inds}^{\infty}
  G_1\left(\frac{x_{\inds} -
      \psi_{\inds}(v)}{v^{\rho_{\tiny \inds}}}\right)G_2\left(\frac{x_{\indt}
      - \psi_{\indt}(v)}{v^{\rho_{{\indt}}}}\right) v^{-2}~\mathrm{d}v,
  \label{eq:link}
\end{equation}
with $\psi_1$ and $\psi_2$ defined by equation \eqref{eq:psi}.\ This
implies that $H(x_{\inds},x_{\indt})=H_1(x_1)\,H_2(x_2)$ if and only
if either $(\kappa_1,\rho_1)=(0,0)$ or $(\kappa_2,\rho_2)=(0,0)$.\
Generally, this condition is not always satisfied and this shows why
random norming is essential for linking conditional independence to
conditioned limit laws.\

Last, we note that conditional independence might not always be an
inherent feature of multivariate extreme value models
\citep{papstro15}.\ Preservation of conditional independence in
conditioned limit laws signifies a potential benefit with this
framework and opens up a possible research direction.

\section{Appendix: Proof of Theorem~\ref{thm:main_thm}}
\label{sec:proof}
The proof of~Theorem~\ref{thm:main_thm} is based on \cite{resnzebe14}
proposition that links kernel convergence to the conditional extreme
value model with deterministic norming.

\begin{lemma}[Proposition 4.1, \cite{resnzebe14}]
  \label{lemma:resnzebe}
  Let $G$ be a non-degererate probability distribution function on
  $[-\infty,\infty)$ and $\pi:(0,\infty)\times
  \mathcal{B}(\mathbb{\bar{R}})\mapsto [0,1]$ be a transition function
  satisfying as $t \rightarrow \infty$
\begin{equation}
  \pi(t,[-\infty,\alpha(t)\,x+\beta(t)]) \overset{\mathcal{D}}{\rightarrow} G(x)
  % \qquad \text{on $[-\infty,\infty]$}.
  \label{eq:kernel_conv_RZ}
\end{equation}
on $[-\infty,\infty]$, where $\alpha(\cdot) \in \mathbb{R}$ and
$\beta(\cdot)>0$ are location and scaling functions.\ There exists a
family of non-degenerate probability distribution functions
$\{G_v:0<v<\infty\}$ on $[-\infty,\infty)$ such that for $0<v<\infty$
and as $t\rightarrow \infty$
\[
K(tv,[-\infty,\alpha(t)\,x + \beta(t)])
\overset{\mathcal{D}}{\rightarrow} G_v(x)
\]
on $[-\infty,\infty]$, if and only if $\alpha,\beta$ are extended
regularly varying with parameters $\rho,\kappa\in\mathbb{R}$,
\[
\frac{\alpha(t x)}{\alpha(t)} \rightarrow x^{\rho} \qquad \text{and}
\qquad \frac{\beta(tx) - \beta(t)}{\alpha(t)} \rightarrow \psi (x)=
\begin{cases}
  \kappa \left(x^{\rho}-1\right)/\rho & \rho \neq 0\\
  \kappa \log x &\rho=0
\end{cases},
\]
for $x>0$.\ In this case, $G_{\inds}=G$, and
\[
\pi(tv_t, [-\infty,\alpha(t)\,x + \beta(t)])\rightarrow
G(\{x-\psi(v)\}/v^{\rho}),
\]
at continuity points of $G$, whenever $v_t=v(t)\rightarrow v
\in(0,\infty)$.
\end{lemma}

\noindent \textbf{Proof of Theorem~\ref{thm:main_thm}}
  The left hand side of expression~\eqref{eq:week1} is equal to
  \begin{IEEEeqnarray}{rCl}
    \lefteqn{ \frac{1}{\mathbb{P}\left(\Xf>t\right)}
      \int_{t}^{\infty}\int_{\mathbb{R}\times \mathbb{R}}
      \,f_{\inds}\left(\frac{x_{\inds}-\locs(t)}{\scas(t)}\right)\,
      f_{\indt}\left(\frac{x_{\indt}-\loct(t)}{\scat(t)}\right)\,
      \mathbb{P}\left(\Xf\in dx_{\scriptsize\indf}, \Xs \in
        dx_{\inds}, \Xt \in
        dx_{\indt}\right)}% \\\\\\ &=& \int_{\inds}^{\infty} \left[
    % \int_{\mathbb{R}} f\left(\frac{y-\beta(t)}{\alpha(t)}\right)
    % \left[\int_{ \mathbb{R}} g\left(\frac{z-b(t)}{a(t)}\right)
    %   \pi_{X,Z}\left( t v,dz \right) \right]
    % \pi_{X,Y}\left( t v,dy \right) \right]
    % \frac{\mathbb{P}\left(X\in
    %   t dv\right)} {\mathbb{P}\left(X>t\right)}
    \nonumber\\\nonumber\\ &=& \int_{\inds}^{\infty} \left[
      \int_{\mathbb{R}}\, f_{\inds}\left(x_{\inds}\right) \left[\int_{
          \mathbb{R}} \,f_{\indt}\left(x_{\indt}\right) \pi_{\indt}\left( t
          v,\scat(t)\,dx_{\indt}+\loct(t) \right) \right]
      \pi_{\inds}\left( t v,\scas(t)\, dx_{\inds} + \locs(t)\right)
    \right] \,\frac{\mathbb{P}\left(\Xf\in t\, dv\right)}
    {\mathbb{P}\left(\Xf>t\right)}.\nonumber
  \end{IEEEeqnarray}
  The term in the inner square brackets of the last expression is
  bounded and, by assumption~\eqref{eq:kernel_conv} and
  Lemma~\ref{lemma:resnzebe}, it converges to
  \[
  \mathbb{E}_{G_{v;\indt}}\{f_{\indt}(\Xt)\} = \int_{\mathbb{R}} f_{\indt}(x_{\indt})
  \,G_{v;\indt}(dx_{\indt}),
  \]
  where $\{G_{v;\indt}:0<v<\infty\}$ is a family of non-degenerate
  probability distributions on $[-\infty,\infty)$ such that
  $G_{1;\indt}=G_{\indt}$. Working similarly for the term in the outer
  square brackets, we get that the left hand side of
  expression~\eqref{eq:week1} equals
  \[
  \int_{\inds}^{\infty} \mathbb{E}_ {G_{v;\inds}}\{f_{\inds}(\Xs)\}\,
  \mathbb{E}_{G_{v;\indt}}\{f_{\indt}(\Xt)\}\, v^{-2} ~dv \qquad
  \text{as $t\rightarrow \infty $}.
  \]
  Expression~\eqref{eq:week2} is obtained in a similar manner and this
  completes the proof.

% Given $\rho,\kappa \in \mathbb{R}$, define the generalized tail
% kernel associated with a distribution $G$ on $[-\infty,\infty]$ as
% the transition function
% $\kappa_G:(0,\infty)\times\mathcal{B}[-\infty,\infty] \mapsto[0,1]$
% given by
% \begin{equation}
%   \kappa_G(y,A) = G(y^{-\rho}[A-\psi(y)]),
%   \label{eq:4.3}
% \end{equation}
% were $\psi$

\subsection*{Acknowledgments}
Ioannis Papastathopoulos acknowledges funding from the SuSTaIn program
- Engineering and Physical Sciences Research Council grant
EP/D063485/1 - at the School of Mathematics, University of Bristol.\
The author would like to thank Sidney Resnick, Kirstin Strokorb and
Jonathan Tawn for helpful discussions. Thanks to Cornell University,
Universit\"{a}t Mannheim and Lancaster University for funding during
visits in 2015.% funding and

\end{document}